\documentclass[10pt,twoside,reqno]{amsart}

\newcommand{\vs}{\vspace}
\newcommand{\bc}{\begin{center}}
\newcommand{\ec}{\end{center}}

\usepackage{amsmath,cite}
\usepackage{amssymb}
\usepackage{amsfonts}
\usepackage{amsthm}

\numberwithin{equation}{section}

\newcommand{\ve}{\varepsilon}
\newcommand{\NN}{\mathbb N^+}
\newcommand{\RR}{\mathbb R}

\begin{document}
\setcounter{page}{1}

\title[Periodic point theorem for generalized graphic contractions]{\large Periodic point theorem for generalized graphic contractions}
\author[Evgeniy Petrov]{Evgeniy Petrov}
\date{}
\maketitle

\vs*{-0.5cm}

\bc
{\footnotesize
Function Theory Department, Institute of Applied Mathematics and Mechanics of the NAS of Ukraine, Batiuka Str. 19, Slovyansk 84116, Ukraine.\\
E-mail: eugeniy.petrov@gmail.com}
\ec

\bigskip

{\footnotesize
\noindent
{\bf Abstract.}
Let $(X,d)$ be a nonempty metric space and let $n\in \NN$. We shall say that $T\colon X\to X$ is a \emph{graphic contraction of order $n$} if there exists $\alpha\in (0,1)$ such that the inequality
$$
  d(T^n x,T^{2n}x) \leqslant \alpha d(x,T^nx)
$$
holds for all $x\in X$. In the case $n=1$ these mapping are known as graphic contractions and are well studied.
In the present paper we establish a theorem on the existence of periodic points for graphic contraction of order $n$. Examples of such mappings having different properties are constructed.



}

\bigskip

\section{Introduction}
Except Banach contraction principle in the Fixed Point Theory on general metric spaces there are several classical fixed point theorems against which metric extensions are usually checked. Perhaps, the best known of them are Nadler's set-valued extension of Banach's theorem, the extension of Banach's theorem to nonexpansive mappings,  Caristi's, Kannan's, Chatterjea's, Suzuki's  theorems, etc. Such generalizations are very numerous and as a rule establish the existence and uniqueness of a fixed point. There are considerably fewer theorems which establish the existence of periodic points of mappings in general metric spaces. The most known is Edelstein's~\cite{Ed62} theorem stating that any $\ve$-contractive mapping of a nonempty compact metric space into itself has a periodic point. See also some generalizations of this theorem~\cite{Ra06,CJU08,Ba66,Se72}. These theorems do not generalize Banach's theorem and have a completely different proof scheme. A theorem on the existence of periodic points for mappings contracting the total pairwise distances between $n$ points was proved in~\cite{P25}. The classical Banach fixed-point theorem was obtained like a simple corollary as well as the fixed-point theorem for mappings contracting perimeters of triangles~\cite{P23}. One more periodic point theorem formulated in terms of $\varepsilon$-$\delta$-conditions was recently obtained by Pant and Rako\v{c}evi\v{c} in~\cite{PR24}.
At the same time, we note that the study of periodic points play an important role in the theory of dynamical systems~\cite{De22} and is also well developed for spaces of special types~\cite{Sa64,Bo71,BGMY80,Fr92,AM65,SY17}. The present paper is a new contribution to the theory of periodic points in general metric spaces. 

Let $(X,d)$ be a metric space and let $T\colon X\to X$ be a mapping such that
\begin{equation}\label{e0}
d(Tx,T^2x)\leqslant \alpha d(x,Tx)
\end{equation}
for some $\alpha\in (0,1)$ and every $x\in X$. In 1968 Rheinboldt~\cite{R68} named mappings satisfying the above inequality~(\ref{e0}) ``iterated contractions'' and established a fixed point theorem for these mappings on a closed subset of $\RR^n$.
In 1972, Rus proved in~\cite{Rus72} the following result.

\bigskip

\noindent{\bf Theorem 1.1.} {\it  
Let $T$ be a continuous self-map of a complete metric space $X$ satisfying condition~\emph{(\ref{e0})}. Then $T$ has a fixed point.
\ }

\bigskip

In 1974, Subrahmanyam~\cite{Su74} established a fixed point theorem for these mappings, which he called ``Banach operators of type k'', considering them on closed subsets of Banach spaces.
In 1979 Hicks and Rhoades~\cite{HR79} proved the existence of a fixed point $x^*$  for such mappings in a complete metric space in the case if the functional $G(X)=d(x,f(x))$ is $T$-orbitally lower semi-continuous at $x^*$. In this connection sometimes such mappings in the literature are called Rus-Hicks-Rhoades mappings~\cite{Park25,park2023relatives,park2023almost,park2024realm}.

Nevertheless, the most widely used name for mappings satisfying inequality~(\ref{e0}) is \emph{graphic contraction}, see~\cite[P.~530]{BPR23}. The class of graphic contractions contains many mappings, including Banach, Kannan, \'{C}iri\'{c}-Reich-Rus, Chatterjea, Zamfirescu, Hardy-Rogers, Berinde, Suzuki contractions, see~\cite[P.~418]{PR19} or~\cite[Sections 5-7]{park2023almost}. See also~\cite{BP23,CPP20,MP18,Fi25,Rus16} for further results related to graphic contractions and their generalizations.

We generalize the concept of graphic contractions as follows.

\bigskip

\noindent{\bf Definition 1.2.} Let $(X,d)$ be a nonempty metric space and let $n\in \NN$. We shall say that $T\colon X\to X$ is a \emph{graphic contraction of order $n$} if there exists $\alpha\in (0,1)$ such that the inequality
  \begin{equation}\label{e1}
   d(T^n x,T^{2n}x) \leqslant \alpha d(x,T^nx)
  \end{equation}
  holds for all $x\in X$.

\bigskip

The following Lemma, which we need below, is a part of the standard proof of Banach contraction principle, in particular its proof follows from the corresponding lemma formulated for $b$-metric spaces,  see~\cite[Lemma 2.2]{MM17} or from the corresponding lemma formulated for semimetric spaces with triangle functions, see~\cite[Lemma 1.2]{PSB24}.

\bigskip

\noindent{\bf Lemma 1.3.} {\it  
Let $(X,d)$ be a metric space. Let $(x_k)_{k\in \NN}$ be a sequence of elements from $X$ having the property that there exists $\gamma\in [0,1)$ such that
\begin{equation}\label{e2}
d(x_{k+1},x_k)\leqslant \gamma d(x_k,x_{k-1})
\end{equation}
for all $k\geqslant 2$. Then the sequence $(x_k)_{k\in \NN}$ is Cauchy.
\ }

\bigskip

\section{The main result}

Let $T$ be a mapping on the metric space $X$. A point $x\in X$ is called a \emph{periodic point of period $n$} if $T^n(x) = x$. The least positive integer $n$ for which $T^n(x) = x$ is called the prime period of $x$, see, e.g.,~\cite[p.~18]{De22}.
By $O_T(x)$ we denote the orbit of the point $x$ under the mapping $T$, $O_T(x)=\{x, Tx, T^2x,\ldots\}$.
Recall that for a given metric space $X$, a point $x \in X$ is said to be an \emph{accumulation point} of  $X$ if every open ball centered at $x$ contains infinitely many points of $X$.

\bigskip

\noindent{\bf Theorem 2.1.} {\it 
Let $(X,d)$ be a complete metric space and $T\colon X\to X$ be a continuous graphic contraction of order $n\in \NN$. Then $T$ has a periodic point. The prime period of this periodic point is a divisor of $n$.
\ }

\bigskip

\noindent{\bf Proof.} If $n=1$, then according to Theorem~1.1 the mapping $T$ has a fixed point, which by definition is a periodic point of prime period one. Let $n\geqslant 2$ and let $x_1\in X$. Consider the iterative  sequence $x_{k+1}=T^k x_{1}$, $k\geqslant 1$  and $n$ of its subsequences:
\begin{itemize}
  \item [$(s^1_k):$] $x_1$, $x_{n+1}= T^nx_1$, $x_{2n+1}= T^{2n} x_1$,  $x_{3n+1}= T^{3n} x_1$, \ldots;
  \item [$(s^2_k):$] $x_2$, $x_{n+2}= T^{n}x_2$, $x_{2n+2}= T^{2n} x_2$,  $x_{3n+2}= T^{3n} x_2$, \ldots;
  \item [$\ldots$]
  \item [$(s^n_k):$] $x_n$, $x_{2n}= T^{n}x_n$, \, $x_{3n}= T^{2n} x_n$, \, \,  $x_{4n}= T^{3n} x_n$, \ldots.
\end{itemize}

In other words $s^i_k=x_{n(k-1)+i}$.
Comparing conditions~(\ref{e1}) and~(\ref{e2}) we get from Lemma~1.3 that every sequence $(s^i_k)$, $i=1,\ldots,n$, is Cauchy. Hence, by the completeness of $X$ every sequence $(s^i_k)$ converges to some $x_i^*\in X$, $i=1,\ldots,n$.

By the continuity of $T$ we have the following implications
$$
\lim\limits_{k\to \infty} s^i_k = x_i^* \,\,  \Rightarrow  \, \,
\lim\limits_{k\to \infty} Ts^i_k = Tx_i^*
$$
for every $i\in 1,\ldots,n$. Observing that $Ts^i_k=s^{i+1}_k$ when $i=1,\ldots,n-1$ and  $Ts^n_k=s^{1}_{k+1}$, we get that
\begin{equation}\label{e6}
Tx_i^* = \lim\limits_{k\to \infty} Ts^i_k = \lim\limits_{k\to \infty} s^{i+1}_k = x_{i+1}^*, \, i=1,\ldots,n-1.
\end{equation}
and
\begin{equation}\label{e5}
Tx_n^* = \lim\limits_{k\to \infty} Ts^n_k = \lim\limits_{k\to \infty} s^{1}_{k+1} = x_{1}^*.
\end{equation}

\textbf{Case a)}. Suppose first that the points $x^*_1, x^*_2,\ldots,x^*_n$ are pairwise distinct. Clearly, in this case $x^*_1, x^*_2,\ldots,x^*_n$ are periodic points of the mapping $T$ of prime period $n$.

\textbf{Case b)}. Suppose that $x^*_1=x^*_2\cdots=x^*_n=x^*$. Clearly, in this case $x^*$ is a fixed point.

We consider the following two cases \textbf{c)} and \textbf{d)} for $n\geqslant 3$, since in the case $n=2$ cases \textbf{c)} and \textbf{d)} are impossible. In case \textbf{c)} we assume that there are equal consecutive points $x_i^*=x^*_{i+1}$, and in case \textbf{d)} we consider the existence of equal points $x_i^*=x^*_{i+p}$ that are not consecutive.

\textbf{Case c)}. Suppose now that for some $i=1,\ldots,n-2$ we have
\begin{equation}\label{e7}
x_i^*=x_{i+1}^*\neq x^*_{i+2}.
\end{equation}
By~(\ref{e6}) we have  $Tx_i^*=x^*_{i+1}$. By~(\ref{e7}) this means that $x_i^*$ is a fixed point of $T$, $Tx_i^*=x_i^*=Tx_{i+1}^*=x_{i+1}^*$. Further, by~(\ref{e6}) we have  $Tx_{i+1}^*=x^*_{i+2}$, which is a contradiction. Hence, case \textbf{c)} is impossible. The cases when $i=n-1$ ($x_{n-1}^*=x_{n}^*\neq x_{1}$), $i=n$ ($x_{n}^*=x_{1}^*\neq x_{2}$) are impossible analogously.

\textbf{Case d)}.
Suppose now that in these set $\{x^*_1, x^*_2,\ldots,x^*_n$\} there are points $x_i^*$ and $x^*_{i+p}$, $p\geqslant 2$, such that $x_i^*=x^*_{i+p}$ and
the points $x_i^*, x^*_{i+1}, \ldots, x^*_{i+p-1}$  are pairwise distinct.

For simplicity, without loss of generality, we can assume that $i=1$. Thus, we consider that
\begin{equation}\label{e8}
x^*_1=x^*_{p+1}
\end{equation}
and the points
$x_1^*, x^*_{2}, \ldots, x^*_{p}$
are pairwise distinct.

By~(\ref{e6}) we have that $Tx^*_k=x^*_{k+1}$ for $k=1,\ldots,p$. Since $x^*_k$ and $x^*_{k+1}$, $k=1,\ldots,p-1$ are different, using equality~(\ref{e8}), we see that $T^p x^*_1=x^*_{1}$. This means that $x^*_1$ is a periodic point of prime period $p$.

Consider the following finite ordered sequences
$$
S_1=(x_1^*, x^*_{2}, \ldots, x^*_{p}),
$$
$$
S_2=(x_{p+1}^*, x^*_{p+2}, \ldots, x^*_{2p}),
$$
$$
S_3=(x_{2p+1}^*, x^*_{2p+2}, \ldots, x^*_{3p}),
$$
$$
\cdots
$$
$$
S_m=(x_{(m-1)p+1}^*, x^*_{(m-1)p+2}, \ldots, x^*_{mp}),
$$
where $m\in \NN$ is the biggest natural number such that
$mp\leqslant n$. Observe that conditions~(\ref{e8}) and~(\ref{e6}) imply that $S_1=S_2=\cdots = S_m$.

We claim that $mp=n$. Indeed, suppose that $mp<n$. Then
$$
1\leqslant |\{x^*_{mp+1}, x^*_{mp+2},\ldots, x^*_{n}\}|=p'<p,
$$
where $|\cdot|$ denotes the cardinality of the set.
By~(\ref{e6}) and~(\ref{e8}) we have
$$
Tx^*_{mp}=x^*_{mp+1}=x^*_1, \quad
x^*_{mp+2}=x^*_2, \quad \ldots, \quad
x^*_{n}=x^*_{mp+p'}=x^*_{p'},
$$
and $Tx^*_{n}=x^*_{p'+1}$. Clearly, $2\leqslant p'+1 \leqslant p$. But by~(\ref{e5}) $Tx_n^*=x_1^*$, which is a contradiction, since $x_1^*, x^*_{2}, \ldots, x^*_{p}$ are pairwise distinct. Thus, $mp=n$ and $p$ is a divisor of $n$.

Clearly, in the cases \textbf{a)} and \textbf{b)} the prime period of the periodic point is also a divisor of $n$ It suffices to note that any combination of consecutive coincidences and dissimilarities between the points $x^*_1, x^*_2,\ldots,x^*_n$ is described by one of the cases \textbf{a)}, \textbf{b)}, \textbf{c)} or \textbf{d)}, which completes the proof.

\bigskip

Theorem~2.1 establishes the existence of at least one periodic point $x$. If the prime period $p$ of this point is strictly greater than 1, then it is clear that the number of periodic points is at least  $p$: all the points of the orbit $O_T(x)$. In the following example, we construct a finite metric space $X$ and a graphic contraction $T\colon X\to X$ of order $6$, which has periodic points with non-intersecting orbits.

\bigskip

\noindent{\bf Example 2.2.} 
Let $(X,d)$ be a metric space such that $X=\{x_1,x_2,x_3,x_4,x_5\}$, $d(x_i,x_j)=1$ for all $i\neq j$ and let a mapping $T\colon X\to X$ be such that: $T(x_1)=x_2$, $T(x_2)=x_1$, $T(x_3)=x_4$, $T(x_4)=x_5$, $T(x_5)=x_3$. It is clear that $x_1$ and $x_2$ are periodic points of prime period $2$; $x_3$, $x_4$ and $x_5$ are periodic points of prime period $3$ and $O_T(x_1)\cap O_T(x_3)=\varnothing$. It is not hard to see that $T$ is a graphic contraction of order $6$: $d(x,T^6x)=d(T^6x,T^{12}x)=0$ for all $x\in X$. Since the metric space $X$ is discrete $T$ is continuous.

\bigskip

\noindent{\bf Example 2.3.} 
Let $a<b$, $a,b\in \RR$ and let
$$
x_n=
  \begin{cases}
        a-\frac{1}{2^n}, & n\text{ is odd;} \\
        b+\frac{1}{2^n}, & n\text{ is even.}
      \end{cases}
$$
for every $n\in \NN$. Let $(X,d)$ be a metric space with $X=\{a,b,x_1,\ldots,x_n,\ldots\}$ and $d$ be the usual Euclidean distance. Define a mapping $T\colon X\to X$ as follows

\begin{equation}\label{e4}
T(x)=
  \begin{cases}
        x_{n+1}, &\text{if } x=x_n, \, \, n\geqslant 1; \\
        b, &\text{if } x=a; \\
        a, &\text{if } x=b.
      \end{cases}
\end{equation}
It is clear that $T$ is not a graphic contraction of order $1$, since
$$
\lim\limits_{n\to\infty}\frac{d(Tx_n,T^2x_n)}{d(x_n,Tx_n)}=1.
$$
Let us show that $T$ is a graphic contraction of order $2$. Consider inequality~(\ref{e1}) with $n=2$ and $x=x_{2k+1}$, $k \in \mathbb N$:
$$
d(x_{2k+3},x_{2k+5})\leqslant \alpha d(x_{2k+1},x_{2k+3}),
$$
$$
\left|\frac{1}{2^{2k+3}}-\frac{1}{2^{2k+5}}\right|\leqslant \alpha \left|\frac{1}{2^{2k+1}}-\frac{1}{2^{2k+3}}\right|,
$$
$$
\frac{24\cdot 2^{2k}}{2^{4k+8}}\leqslant \alpha \frac{6\cdot 2^{2k}}{2^{4k+4}},
$$
$$
\frac{1}{4}\leqslant \alpha.
$$
Thus, in this case inequality~(\ref{e1}) holds for $\alpha \in [\frac{1}{4},1)$. Analogously,  inequality~(\ref{e1}) holds for $\alpha \in [\frac{1}{4},1)$, $n=2$ and $x=x_{2k}$, $k \in \NN$.
Clearly, for $x=a$ or $x=b$ inequality~(\ref{e1}) holds for any $\alpha \in [0,1)$. Thus, $T$ is a graphic contraction of order $2$ for any $\alpha \in [\frac{1}{4},1)$ with two periodic points $a$ and $b$ of prime period $2$.

Clearly, the mapping $T$ is continuous at all isolated points of the space $X$. Note that $X$ contains only two accumulation points, $a$ and $b$. Observe that if $\xi_k \to a$, then $T\xi_k \to Ta = b$, and if $\xi_k \to b$, then $T\xi_k \to Tb = a$. Therefore, $T$ is continuous at every point in the space $X$.

\bigskip

Next we consider a modification of the previous example to illustrate case d) considered in the proof of Theorem~2.1.

\bigskip

\noindent{\bf Example 2.4.} 
Let $a<b$, $a,b\in \RR$ and let $(x_n)$ be a sequence consisting  of
the following subsequences:
\begin{itemize}
  \item [$(s^1_k):$] ($x_1$, $x_{5}$, $x_{9}$, \ldots) = ($a-\frac{1}{2^1}$, $a-\frac{1}{2^2}$, $a-\frac{1}{2^3}$, \ldots);
  \item [$(s^2_k):$] ($x_2$, $x_{6}$, $x_{10}$, \ldots) =
      ($b+\frac{1}{2^1}$, $b+\frac{1}{2^2}$, $b+\frac{1}{2^3}$, \ldots);
  \item [$(s^3_k):$] ($x_3$, $x_{7}$, $x_{11}$, \ldots) =
      ($a-\frac{1}{3^1}$, $a-\frac{1}{3^2}$, $a-\frac{1}{3^3}$, \ldots);
  \item [$(s^4_k):$] ($x_4$, $x_{8}$, $x_{12}$, \ldots) =
      ($b+\frac{1}{3^1}$, $b+\frac{1}{3^2}$, $b+\frac{1}{3^3}$, \ldots).
\end{itemize}
Let $(X,d)$ be a metric space with $X=\{a,b,x_1,\ldots,x_n,\ldots\}$, where $d$ is the usual Euclidean distance. Define a mapping $T\colon X\to X$ as in~(\ref{e4}).
It is clear that
$$
x^*_1=\lim\limits_{k\to\infty}s^1_k =a, \quad
x^*_2=\lim\limits_{k\to\infty}s^2_k =b, \quad
x^*_3=\lim\limits_{k\to\infty}s^3_k =a, \quad
x^*_4=\lim\limits_{k\to\infty}s^4_k =b.
$$
Similarly to the previous example we can establish that $T$ is a graphic contraction of order $4$ for any $\alpha \in [\frac{1}{2},1)$ with two periodic points $a$ and $b$ of prime period $2$. Analogously, $T$ is not a graphic contraction of order $1$. $T$ is not a graphic contraction of order $3$ since $d(a,T^3a)=d(T^3a,T^6a)=b-a$.

Let us show that $T$ is not a graphic contraction of order $2$. In order to show that inequality~(\ref{e1}) does not hold for some $\alpha\in [0,1)$ and all $x\in X$ with $n=2$ it is sufficient to show that
$$
\lim\limits_{k\to \infty} \frac{d(T^2x_{n_k},T^4x_{n_k})}{d(x_{n_k},T^2x_{n_k})}=1
$$
for some subsequence $(x_{n_k})$ of the sequence $(x_n)$. Let $x_{n_k}=x_{4k-1}=a-\frac{1}{3^k}$, $k\in \NN$.
Note that
$T^2x_{n_k}=T^2x_{4k-1}=a-\frac{1}{2^{k+1}}$, $T^4x_{n_k}=T^4x_{4k-1}=a-\frac{1}{3^{k+1}}$.
Further,
$$
\lim\limits_{k\to \infty}\left|\frac{1}{3^{k+1}}-\frac{1}{2^{k+1}}\right|\colon
\left|\frac{1}{2^{k+1}}-\frac{1}{3^k}\right|
=\lim\limits_{k\to \infty}
\left|\frac{3\cdot 3^k-2\cdot 2^{k}}{6\cdot6^k}\right|\colon
\left|\frac{3^k-2\cdot 2^k}{2\cdot 6^k}\right|
$$
$$
=\lim\limits_{k\to \infty}
\left|\frac{3-2\cdot (2/3)^{k}}{6}\right|\cdot
\left|\frac{2}{1-2\cdot(2/3)^k}\right|=1.
$$
Thus, $T$ is a graphic contraction of order $4$ with periodic points of prime period $2$.

\bigskip

\noindent{\bf Example 2.5.} 
As it was mentioned in the introduction, many known mappings such as Banach, Kannan, Chatterjea, etc., are graphic contractions (of order $1$). Let $(X,d)$ be a metric spaces and $T\colon X\to X$ be a mapping such that $T^n$ is a graphic contraction. Substituting $T^n$ in~(\ref{e0}) instead of $T$ we get inequality~(\ref{e1}). Thus, numerous examples of graphic contractions of order $n$ include mappings whose $n$-th iterations act as Banach~\cite{Br68} $d(T^nx,T^ny)\leqslant \alpha d(x,y)$, Kannan~\cite[Corollary 3.3]{Go18} $d(T^nx,T^ny)\leqslant \alpha (d(x,T^nx)+d(y,T^ny))$, Chatterjea $d(T^nx,T^ny)\leqslant \alpha (d(x,T^ny)+d(y,T^nx))$, or other types of contractions.
Note that according to~\cite[Lemma 3.2]{Go18}, if the $n$-th iteration $T^n$ has only one fixed point, then the mapping $T$ itself also has only one fixed point.

\bigskip

\textbf{Acknowledgement.} E. Petrov was supported by the National Research Foundation of Ukraine, project No. 2025.07/0369 ``Qualitative methods of nonlinear analysis of heterogeneous structures''.

\end{document}